# A Ratio to Evaluate Harvest Procedures Management in an Economic System where Resources Dynamics is ruled by an Ornstein- Uhlenbeck Process

M. A. M. Ferreira[1], J. A. Filipe

Instituto Universitário de Lisboa (ISCTE – IUL), BRU - IUL, Lisboa, Portugal

**Abstract**

The assessing resources dynamics problem, in the context of an economic system with Gaussian consumption and deterministic productivity, is considered in this paper. Basically it is presented a discrete time recursive equation that supports the recourse to the Ornstein-Uhlenbeck diffusion process. Some assumptions on the regeneration of the process are made, in order to observe the system equilibrium in what concerns the resources depreciation or accumulation. The objective of this work is to present a result on the sign of a ratio that can be used to evaluate harvest procedures in this context.

## 1 Introduction

A specific kind of harvesting rules[2] of an economic system, where resources dynamics are described by an Ornstein-Uhlenbeck[3] diffusion process, that is, the diffusion process satisfying the stochastic differential equation

---

[1] manuel.ferreira@iscte.pt

[2] About resources harvesting see for instance [3].

[3] The Ornstein–Uhlenbeck process, named after Leonard Ornstein and George Eugen Uhlenbeck, see [11], is a stochastic process that, roughly speaking, describes the velocity of a massive Brownian particle under the influence of friction. It is stationary, Gaussian and Markovian, being the only nontrivial process that satisfies these three conditions, up to allowing linear transformations of the space and time variables. Over time, the process tends to drift towards its long-term mean: such a process is called mean-reverting. This process can be thought as a modification of the random walk in continuous time, or Wiener process, in which the properties of the process have been changed so that there is a tendency of the walk to move

$$dX(t) = (a + bX(t))dt + dB(t), \quad X(0) = x \quad (1.1),$$

where $B(t)$ is a Wiener process, are studied along this work.

The problem to be studied in this paper can be outlined as follows:

- Supposing that a system like this one is explored according to a particular kind of management or harvesting rules, to which conclusions is it possible to arrive at about the amount of resources made available?

It is assumed the politics that every time the resources process hits the set $\{\eta, \theta\}$, one collects an amount $\eta - x$ or $\theta - x$, with values $Q(\eta)$ or $Q(\theta)$, respectively. This makes the process return instantaneously to the initial level $x$. Also it is made the assumption that the procedure is continuously repeated. Then the economic resources stock is represented effectively through a $X(t)$ modification, called $\tilde{X}(t)$, achieved by regeneration at states $\eta$ and $\theta$.

Call $N(t)$ the regeneration epochs counting process. So the cumulative process

$$R(t) = \sum_{n=0}^{N(t)} Q_n \quad (1.2)$$

represents the accumulated resources harvested in such an economic environment. In fact $\{Q_n\}$ is a sequence of i. i. d. random variables, assuming values on the set $\{Q_n(\eta), Q_n(\theta)\}$, representing the harvesting procedure flow. The main objective is to conclude on the quality of the adopted management scheme[4].

The random variable $T$ characterizes the sequence of the recurrence times of the economic resources process $\tilde{X}(t)$ on the set $\{\eta, \theta\}$. It corresponds to the travel time of $X(t)$ from $x$ to $\{\eta, \theta\}$.

Specially, a result on the sign of the ratio

---

back towards a central location, with a greater attraction when the process is further away from the center. It can also be considered as the continuous-time analogue of the discrete-time AR (1) process.

The Ornstein–Uhlenbeck process is one of several approaches used to model with modifications, interest rates, currency exchange rates, and commodity prices stochastically. In this case

$$dX(t) = \theta(\mu - X(t))dt + \sigma dW(t); \quad \theta, \mu, \sigma > 0.$$

The parameter μ represents the equilibrium or mean value supported by fundamentals; σ the degree of volatility around it caused by shocks, and $\theta$ the rate by which these shocks dissipate and the variable reverts towards the mean. One application of the process is a trading strategy known as "pairs trade".

[4] With the help of the classical renewal theorem, see [9], establishing:

$$\lim_{t \to \infty} \frac{R(t)}{t} = \frac{E[Q]}{E[T]}, \quad a.s.$$

$$\gamma(\theta) = \lim_{\delta \to \Delta} \frac{E[Q]}{E[T]} \qquad (1.3),$$

where $\delta$ and $\Delta$ are variables to be defined later, depending on $\eta, x$ and $\theta$, is highlighted.

In the next section, more results about the process in (1.1) are presented. The main result is shown in the last section. The references [1, 2, 4, 5, 6, 7, 10] are important sources for this theme.

## 2 Detailing the Resources Dynamics Model

Consider the economic system with resources dynamics suggested by the stochastic difference equation

$$Y_0 = x, Y_{nh} = e^{bh}(Y_{(n-1)h} + W_{nh}), \quad n = 1, 2, \ldots \qquad (2.1),$$

where $\{W_{nh}\}$ is a sequence of *i. i. d.* random variables with Gaussian distribution, mean $ah$ and variance $h$; $a < 0$ and $b, h > 0$.

So $\{Y_{nh}\}$ is also a Gaussian sequence with parameters

$$\begin{aligned} E[Y_{nh}] &= xe^{nbh} + ahe^{bh}\frac{1-e^{nbh}}{1-e^{bh}}, \\ Cov[Y_{mh}, Y_{nh}] &= he^{(2+n-m)bh}\frac{1-e^{2mbh}}{1-e^{2bh}}, \quad 0 \leq m \leq n, \end{aligned} \qquad (2.2).$$

Behind (2.1) is the idea of an economic system with a deterministic productivity, multiplying the stock available after Gaussian consumption, in each period of time.

In the identities

$$\begin{aligned} E[Y_{t+h} - Y_t | Y_t = y] &= (a + by)h + o(h), \\ E[(Y_{t+h} - Y_t)^2 | Y_t = y] &= h + o(h), \end{aligned} \qquad (2.3),$$

it is suggested the approximation of $\{Y_{nh}\}$, as $h \to 0$, by a continuous time process satisfying the equation (1.1). They are obtained computing the left side McLaurin's expansion for $t = nh$.

The solution of (1.1) is the Ornstein-Uhlenbeck process

$$X(t) = e^{bt}\left(x + a\int_0^t e^{-bs}ds + \int_0^t e^{-bs}dB(s)\right) \qquad (2.4),$$

with parameters

$$\begin{aligned} E[X(t)] &= xe^{bt} + \frac{a}{b}(e^{bt} - 1), \\ Cov[X(s), X(t)] &= e^{b(s+t)}\frac{1-e^{2bs}}{2b}, \quad 0 \leq s \leq t, \end{aligned} \qquad (2.5).$$

The following functionals concerning the Ornstein-Uhlenbeck process will be needed later in this text:

- Let $0 \leq \eta < x < \theta$ and recall the standard Gaussian distribution function: $\Phi(u) = \int_{-\infty}^{u} \phi(v) dv$ with $\phi(v) = \frac{1}{\sqrt{2\pi}} e^{-\frac{v^2}{2}}$.

- Call $\rho(x)$ the probability that $X(t)$ reaches $\eta$ before $\theta$, coming from $x$. The evaluation of this functional in general diffusion processes is a well-known result, see for instance [8]:

$$\rho(x) = \frac{\Phi(\alpha + \beta\theta) - \Phi(\alpha + \beta x)}{\Phi(\alpha + \beta\theta) - \Phi(\alpha + \beta\eta)}, \quad \alpha = \frac{a\sqrt{2}}{\sqrt{b}} \text{ and } \beta = \sqrt{2b} \quad (2.6).$$

- Calling $\psi(x)$ the expected travel time value, from $x$ to $\{\eta, \theta\}$ [5]:

$$\psi(x) = \frac{2}{\beta}\left(\rho(x)\int_{\eta}^{x}\frac{\Phi(\alpha + \beta u) - \Phi(\alpha + \beta\eta)}{\Phi(\alpha + \beta u)} du \right.$$

$$\left. + (1 - \rho(x))\int_{x}^{\theta}\frac{\Phi(\alpha + \beta\theta) - \Phi(\alpha + \beta u)}{\Phi(\alpha + \beta u)} du\right) \quad (2.7).$$

## 3 The Harvest Procedures Evaluation Ratio

Following the mentioned earlier, the determinant source of information on the managing procedures, observed in this work, is the ratio

$$\frac{E[Q]}{E[T]} = \frac{Q(\theta) + (Q(\eta) - Q(\theta))\rho(x)}{\psi(x)} \quad (3.1).$$

So a result on the sign of

$$\gamma(\theta) = \lim_{\delta \to \Delta} \frac{E[Q]}{E[T]} \quad (3.2),$$

that enhances the way the described economic system is fragile, as to what concerns the chosen harvesting technique, is now derived:

**Proposition 3.1**

A) If $\delta = x, \Delta = \eta$ and $Q(y) = y - x$, the functions

$$\gamma(\theta) \text{ and } \Phi(\alpha + \beta\eta) + \beta \phi(\alpha + \beta\eta)(\theta - \eta) - \Phi(\alpha + \beta\theta)$$

share the same sign. And in particular $\gamma(\theta)$ is positive whenever $-\frac{\alpha}{\beta} \leq \eta < \theta$.

B) If $\delta = x, \Delta = \theta$ and $Q(y) = y - x$, the functions

$$\gamma(\theta) \text{ and } \Phi(\alpha + \beta\theta) + \beta \phi(\alpha + \beta\theta)(\eta - \theta) - \Phi(\alpha + \beta\theta)$$

---

[5] For this functional evaluation in general diffusion processes see for instance again [8].

share the same sign.

**Dem**.: For A) noting that

$$\gamma(\theta) = \frac{-1 - (\theta - \eta)\rho'(\eta)}{\psi'(\eta)} \quad (3.3),$$

the conclusion follows from the fact that $\psi'(\eta)$ is always positive and the numerator of (3.3) has the same sign as the function $\Phi(\alpha + \beta\eta) + \beta\phi(\alpha + \beta\eta)(\theta - \eta) - \Phi(\alpha + \beta\theta)$[6].

Now, for B)

$$\gamma(\theta) = \frac{1 + (\theta - \eta)\rho'(\theta)}{-\psi'(\theta)} \quad (3.4).$$

The conclusion follows from the fact that $\psi'(\theta)$ is always negative and that the numerator of (3.4) has the same sign as the function $\Phi(\alpha + \beta\theta) + \beta\phi(\alpha + \beta\theta)(\eta - \theta) - \Phi(\alpha + \beta\theta)$[7]. ∎

## 4 Concluding Remarks

Proposition 3.1 is a result that allows evaluating the adopted management scheme, referred behind, quality. It is raised with very simple mathematical tools. It may be pointed some weakness on it since the sign of $\gamma(\theta)$ is determined there indirectly. But this does not invalidate its practical utility resulting from its simplicity. Being a difficult task to obtain a practical result to determine directly the $\gamma(\theta)$ sign, it is remarkable that this result relates it with the sign of functions for which it is much simpler to determine.

---

[6] As for this conclusion note that $\Phi(\alpha + \beta\eta) + \beta\phi(\alpha + \beta\eta)(\theta - \eta)$ is the slope of the tangent to the function $\Phi(\alpha + \beta\theta)$ graph at $\eta$.

[7] As for this conclusion note that $\Phi(\alpha + \beta\theta) + \beta\phi(\alpha + \beta\theta)(\eta - \theta)$ is the slope of the tangent to the function $\Phi(\alpha + \beta\eta)$ graph at $\theta$.